\documentclass[11pt,a4paper]{article}

\usepackage[T1]{fontenc}
\usepackage[utf8]{inputenc}
\usepackage{lmodern}
\usepackage{amsmath,amssymb,amsthm}
\usepackage{mathtools}
\usepackage{geometry}
\usepackage{hyperref}
\usepackage{graphicx}
\usepackage{xcolor}
\usepackage{enumitem}
\usepackage{comment}
\usepackage{subcaption}
\usepackage{placeins}

\hypersetup{
	colorlinks=true,
	linkcolor=blue,
	citecolor=blue,
	urlcolor=blue,
	pdftitle={Paper 4},
	pdfauthor={Johannes Hagel},
}

\newtheorem{theorem}{Theorem}[section]

\theoremstyle{definition}
\newtheorem{definition}[theorem]{Definition}
\newtheorem{remark}[theorem]{Remark}

\DeclareMathOperator{\Disc}{Disc}

\title{Algebraic Detection of Tube Rupture via a Cubic Discriminant Criterion}

\author{
	Johannes Hagel\\
	\small Alexander--von--Humboldt--Gymnasium, Neuss, Germany\\
	\small \texttt{johannes.hagel@gmail.com}
}

\date{December 2025}
\begin{document}
\maketitle

\begin{abstract}
	We investigate the rupture of invariant tubes in a class of nonautonomous dynamical systems arising from time-dependent Ermakov-type equations.
	
	Starting from an exactly tube-integrable reference system, we analyze a time-dependent invariant obtained from a positivity-preserving second-order perturbative construction, which provides a near-integrable geometric description of the dynamics. While this approximation does not preserve exact invariance, its algebraic structure remains sufficiently robust to allow for a precise characterization of tube opening and loss of confinement.
	
	For fixed time, the discriminant of the approximate invariant with respect to the   momentum variable defines a cubic polynomial in the configuration variable. We show that the invariant tube admits an unbounded bridge if and only if the associated  cubic possesses exactly one real root. This yields a purely algebraic rupture criterion based on the cubic discriminant and reduces the full geometric problem to the evaluation of a single scalar function of time.
	
	Applying this criterion reveals a sequence of isolated bridge windows whose temporal organization undergoes a transition from one opening per \(2\pi\) cycle to two openings per cycle, corresponding to a period-halving in time. These windows can be represented compactly by a one-dimensional box-plot visualization, which faithfully captures the underlying geometry and highlights the progressive densification and widening of escape-enabling intervals.
	
	The results demonstrate that algebraic diagnostics derived from time-dependent invariants can retain sharp predictive power for rupture phenomena even when exact tube integrability is weakly perturbed.
\end{abstract}

\bigskip \bigskip


\newpage

\section{Introduction}

Time--dependent Hamiltonian systems occupy an intermediate position between classical autonomous integrable dynamics and fully chaotic behavior.
While autonomous integrable systems are characterized by time--independent invariants and a fixed foliation of phase space, explicit time dependence generically destroys exact integrability and allows for a much richer geometric evolution of invariant structures.
Understanding how ordered motion can persist, deform, or eventually break down in such systems remains a central problem in nonlinear dynamics.

Classical approaches to nonautonomous dynamics typically rely on Floquet theory,
averaging, or perturbative KAM-type arguments \cite{Floquet1883,Arnold1989,SandersVerhulstMurdock2007}.

In a series of recent works, a new class of systems was introduced in which the dynamics remains geometrically constrained despite explicit time dependence \cite{Hagel2025Poly,Hagel2025TubeInteg,Hagel2025TubeRupt1}.

The central concept is that of \emph{tube integrability}: instead of invariant tori, the motion is confined to time--dependent invariant tubes in an extended phase space.
These tubes are not static objects; their geometry evolves with time and may undergo topological transitions.
In particular, solutions may remain bounded for long times before escaping through temporary openings of the invariant tube, a phenomenon referred to as \emph{tube rupture}.
Unlike classical separatrix crossing, tube rupture is inherently geometric and cannot be detected by local stability analysis alone.

A remarkable feature of tube--integrable and near
integrable systems is that their behavior differs qualitatively from that of classical autonomous integrable systems.
The invariants governing the motion are explicitly time dependent, and their coefficients may evolve slowly or even secularly.
As a consequence, solutions can be guided adiabatically through a sequence of distinct geometric regimes, giving rise to transient superpositions, intermittent openings, and time--localized structural rearrangements of the invariant tubes.
To our knowledge, such phenomena do not occur in autonomous integrable systems and appear to be intrinsically linked to the explicit time dependence and the slow modulation of invariant coefficients.
We conjecture that this type of behavior is generic for a broad class of nonautonomous tube--integrable systems and represents a structural mechanism distinct from both classical integrability and chaotic dynamics, as suggested by the tube–integrable reference framework introduced in \cite{Hagel2025TubeInteg}.”

The present paper is devoted to a systematic and purely algebraic characterization of tube rupture.
Our starting point is a time--dependent invariant that is quadratic in the momentum variable, as obtained from second--order perturbation theory with strict positivity preservation.
For fixed time, the discriminant of this invariant with respect to the momentum defines a real function $\Disc(z,t)$ that is cubic in the configuration variable $z$.
The existence or absence of a tube opening can then be reduced to a simple algebraic criterion: a bridge allowing unbounded motion exists if and only if the cubic equation $\Disc(z,t)=0$ possesses exactly one real root.
This observation leads naturally to a second discriminant, taken with respect to $z$, which provides a compact and computationally efficient criterion for the detection of bridge windows.

This reduction has several important consequences.
First, it allows for a precise localization of the onset of tube rupture in time.
Second, it reveals a nontrivial temporal organization of bridge events, including a transition from one opening per $2\pi$ period to two openings per $2\pi$, corresponding to a period--halving in time.
Third, the purely algebraic nature of the criterion makes it possible to condense the full geometric information into a one--dimensional time representation, visualized here in the form of a box plot of bridge intervals.
This representation highlights a progressive densification of escape--enabling intervals as time increases.

It is important to emphasize that the present analysis is geometric and algebraic in nature.
No probabilistic assumptions are made, and no claims regarding Lyapunov stability are required.
In particular, while the $z$--dynamics is constrained by tube--integrability, no conserved quantity is known for the auxiliary $y$--equation governing the time dependence of the invariant coefficients.
The possibility that regular tube--confined motion in $z$ may coexist with irregular or weakly chaotic behavior in the driving subsystem therefore remains open and will be addressed in future work.

The paper is organized as follows.
In Section~2, we recall the construction of the time-dependent invariant underlying tube integrability and introduce the perturbative approximation used throughout the analysis.
Section~3 develops the algebraic discriminant formalism and derives the cubic criterion for the existence of bridge openings.
The temporal organization of bridge windows, including the transition from one opening per \(2\pi\) cycle to a period-halved structure, is analyzed and visualized using a compact box-plot representation.
Finally, Section~4 summarizes the main results and discusses their implications for near integrable, time-dependent dynamics.

\section{Example of a nonautonomous integrable system}

Given the second order, nonlinear and nonautonomous ODE
	\begin{equation}
		z''(t) + \omega^2 z(t) +  g(t)\,z(t)^2 = 0,
		\label{eq:tube-physical}
	\end{equation}
	we identified an exactly integrable case of this equation with an invariant $I(z,p,t) \, \, , \, \,  p=z'$, if the following four conditions (S1)-(S4) hold \cite{HagelBouquet1992}:

		\begin{align*}
			(S1)\;& z'' + \omega^2 z + g(t) z^2 = 0 \\[0.5em]
			(S2)\;& g(t) = \alpha_2(t)^{-5/2} \\[0.5em]
			(S3)\;& \alpha_2''' + 4\omega^2 \alpha_2' - [C_1\cos(\omega t)+C_2\sin(\omega t)]\,\alpha_2^{-5/2} = 0 \\[0.5em]
			(S4)\;& I(z,p,t) = \alpha_2(t)p^2 - \alpha_2'(t)zp + \alpha_1(t)p \\
			&\quad + \Big(\omega^2\alpha_2(t)+\tfrac{1}{2}\alpha_2''(t)\Big)z^2
			- \alpha_1'(t) z + \tfrac{2}{3}\alpha_2(t) g(t) z^3 \;=\; K
		\end{align*}
	
	\begin{remark}[Historical note on the Bouquet--Hagel system]
		The nonautonomous system defined by conditions (S1)--(S4) was developed jointly by
		C.~Bouquet and J.~Hagel, \cite{HagelBouquet1992} in the mid--1980s in the context of early studies on nonlinear
		time--dependent oscillators.
		The construction emerged from discussions during a conference on nonlinear dynamics
		in 1986 and was subsequently documented in a joint internal report prepared at CERN.
		
		Although the system was not published in a refereed journal at that time, the report
		contains the original derivation of the invariant structure and the associated auxiliary
		equation governing the time dependence of the coefficients.
		In the present work, we therefore refer to this exactly integrable reference case as the
		\emph{Bouquet--Hagel system}, acknowledging its joint origin and historical development. 
		The present attribution is based on the available documentation and correspondence;
		any earlier or independent derivations, if identified, will be acknowledged in future work.

	\end{remark}

	Letting $C_2$=0 and $C_1=\varepsilon$ and renaming
	\begin{equation}
		y(t) = \alpha_2(t) \qquad , \qquad
	\end{equation}
	Equation \emph{(S3)} becomes
	\begin{equation}
		y''' + 4 \omega^2\,y'
		- \varepsilon \cos \omega t \, y^{-5/2} = 0.
	\end{equation}
	The case $C_2=0$ means no loss of generality since every superposition of the form $C_1\cos{\omega t}+C_2 \sin{\omega t}$ can be written as a single cosine function containing a certain phase:, hence $C \cos{(\omega t+\Phi)}$ with  $C=\sqrt{C_1^2+C_2^2}$.

	\subsection{Tube integrability}
	
	We consider a nonautonomous second--order equation
	\begin{equation}
		z'' = F(z,p,t), \qquad p = z',
	\end{equation}
	and assume the existence of a smooth real function
	\begin{equation}
		I(z,p,t)
	\end{equation}
	such that every solution $(z(t),p(t))$ satisfies
	\begin{equation}
		I(z(t),p(t),t) = K,
	\end{equation}
	for a constant $K$ determined by the initial data.
	
	\begin{definition}[Tube integrability]
		The system is called \emph{tube integrable} if the following two conditions hold:
		\begin{enumerate}
			\item \textbf{Invariant surface:}  
			A smooth invariant relation $I(z,p,t)=K$ exists for all solutions.
			
			\item \textbf{Aperiodic time dependence:}  
			The coefficient functions in $F(z,p,t)$, and therefore the invariant 
			$I(z,p,t)$ itself, are \emph{not} periodic in~$t$.
		\end{enumerate}
		In this case the invariant set
		\begin{equation}
			\mathcal{T}_K = \{(z,p,t)\in \mathbb{R}^3 \;|\; I(z,p,t)=K\}
		\end{equation}
		forms a smooth, nonclosing two--dimensional surface in the extended phase
		space $(z,p,t)$, topologically equivalent to an infinite cylinder.  
		We refer to~$\mathcal{T}_K$ as an \emph{invariant tube} \cite{Hagel2025TubeInteg}.
	\end{definition}

	\begin{remark}
		If the time dependence were periodic, the invariant surface $I=K$ would close
		onto itself after one period, giving rise to a torus geometry.  
		Tube integrability therefore represents the genuinely nonautonomous
		counterpart of classical torus--based integrability: the system possesses an
		invariant, but aperiodicity in~$t$ forces the solution to evolve along an
		unclosed invariant tube rather than a torus.
	\end{remark}
	\begin{remark}[No boundedness requirement]
		Tube integrability, as defined above, does not require the solutions 
		$(z(t),p(t))$ to remain bounded for all~$t$.  
		The invariant relation $I(z,p,t)=K$ may give rise to invariant tubes 
		$\mathcal{T}_K$ that are noncompact or even unbounded in the $(z,p)$--directions.  
		Such behaviour does not violate integrability: the identity 
		\[
		I(z(t),p(t),t)=K
		\]
		continues to hold, and the geometric structure of an invariant tube remains 
		well defined.
		
		In particular, the algebraic mechanism of \emph{rupture} --- detected through 
		changes of sign in the discriminant of $\partial I/\partial p$ --- may open or 
		split an invariant tube without destroying integrability.  
		Unlike the breakdown of invariant tori in autonomous Hamiltonian systems, 
		rupture does not signal the loss of an integral; it merely alters the 
		topology of the invariant surface within the extended phase space.  
		Thus, unbounded or opening tubes are fully compatible with tube integrability. This mechanism is fundamentally different from classical separatrix crossing
		in autonomous Hamiltonian systems \cite{Arnold1989}.
		
	\end{remark}
	\begin{remark}[Integrability versus tube integrability]
		Classical (Liouville) integrability implies the existence of a full set of
		first integrals in involution and, after a canonical transformation,
		reduces the dynamics to quadratures on an invariant torus.
		In this sense, ``integrable'' is often interpreted as ``explicitly solvable''
		or ``separable''.
		
		Tube integrability is conceptually different. 
		The existence of a smooth invariant 
		\[
		I(z,p,t)=\mathrm{const}
		\]
		reduces the second--order equation to a nonautonomous 
		first--order ODE, but the reduced equation is, in general,
		\emph{not} separable. 
		Hence tube integrability guarantees a geometric
		integrability (motion on a two–dimensional invariant tube in $(z,p,t)$),
		but it does not imply closed-form solutions nor classical separability.
	\end{remark}
	In the following sections we will not further exploit action–angle constructions. Instead, we focus on the behavior of algebraic invariants obtained from perturbative approximations of the exact tube–integrable structure described above.

 	
\section{Analytic computation of the rupture thresholds}
If an invariant $I(z,p,t)=K$ exists in a non autonoumous system (Hamiltonian or other), for certain values of the initial 
conditions there might exist the case of a limited region of closed curves representing bounded motion depending on the initial conditions and on the independent variable (time).
Exceeding these limits then leads to unbounded motion. However these limits do not mean loss of integrability as the exact solution
still moves on the invariant surface. The only difference lies in the topology of the invariants. Hence, in case of integrability the question of boundedness of the solution $z(t)$ can be answered by a purely algebraic investigation of the invariant function $I$

\paragraph{Convention (\texorpdfstring{$\omega=1$}{omega=1}).}
Throughout this paper we work in the physical time variable $t$ and restrict to the
normalized frequency $\omega=1$ in the $z$--equation. This choice is made solely for
notational simplicity and keeps all formulas in their shortest form. The general case
$\omega\neq 1$ is obtained by the standard rescaling $\tau=\omega t$ (with the corresponding
parameter redefinitions), but is not pursued here.

Our starting point is the perturbative solution of the $y$–system
given in Eqs.~(S1)–(S4) in \cite{Hagel2025TubeInteg}.
. We truncate $y(t)$ at order
$\mathcal{O}(\varepsilon^{2})$,
\begin{equation}
y(t) \;\approx\; y_0 \exp\!\bigl(
\varepsilon \rho_1(t) + \varepsilon^2 \rho_2(t)+O(\varepsilon^3)
\bigr) \label{eq:ypert}
\end{equation}
where 
\begin{equation}
	\rho_1(t)=y_0^{-7/2} \left(\frac{1}{3}\sin{t}-\frac{1}{6}\sin{2t} \right) \label{eq:rho1}
\end{equation}
\begin{equation}
	\rho_2(t)=y_0^{-7/2} \left( -\frac{5}{288}-\frac{1}{24}\cos{t}+\frac{19}{288}\cos{2t}-\frac{1}{72}\cos{3t} \\
	+\frac{1}{144}\cos{4t}+\frac{5}{96} t \sin{2t} \right)
\end{equation}
and use this expression to construct the time–dependent coefficient
$g(t)=y(t)^{-5/2}$ entering the $z$–equation. 

It is very important at this point to mention two important facts:

\vspace*{0.5cm}

\noindent \underline{\textbf{First:}}

\noindent Truncating the perturbation expansion of the equation $y'''+4y'=\varepsilon y^{-5/2} \cos{t}$ at any finite order means that the invariant defined in Eqs. (S1)-(S4) will no longer remain exactly constant as it is in the associated Bouquet-Hagel system. In this case our system, due to truncation of the perturbation series becomes weakly nonintegrable. Despite of this fact, we will still use the near invariant to determine the approximate geometric behaviour of the associated z-equation given by Eq. (S1) \cite{Hagel2025TubeRupt1}.

\vspace*{0.5cm}

\noindent \underline{\textbf{Second:}}

It must be noted that the exponential form of the perturbation series as given in (\ref{eq:ypert}) is essential, since disregarding  truncation of the series in the exponent, this method exactly preserves the property $y(t)>0$ which prevents $y^{-5/2}$ from becoming imaginary. Hence expanding $\log(y(t))$ keeps $y(t)$ a physical meaningful quantity.

Figures~1--2 clearly demonstrates this fact by comparing the numerical solution of the $y$-equation with two different perturbative approximations for the standard parameter set $\varepsilon=0.08$ and $y_0=1$.
Figure~1 shows the reference solution obtained by direct numerical integration of the full equation over a long time interval.
This solution remains strictly positive and exhibits a slow secular modulation superimposed on the dominant oscillatory behaviour.

In Figure~2, the perturbative solution is represented in exponential form,
\[
y(t)=\exp\!\big(\ln y_0 + \varepsilon \rho_1(t) + \varepsilon^2 \rho_2(t)\big),
\]
which preserves positivity by construction.
Remarkably, this representation reproduces both the amplitude and the phase of the numerical solution over the entire integration interval with high accuracy.
In particular, no spurious sign changes or artificial amplitude drifts are observed.

Figure~2 displays the straightforward perturbative expansion obtained by a direct Taylor series in $\varepsilon$ up to second order.
Although this approximation performs reasonably well at early times, it gradually loses accuracy as time increases.
In contrast to the exponential formulation, positivity of $y(t)$ is not guaranteed, and noticeable deviations in amplitude and phase accumulate.
This comparison clearly demonstrates that the exponential perturbative form provides a substantially more robust global approximation than the naive Taylor expansion.

\vspace*{0.5cm}

\begin{figure}[!htbp]
	\centering
	\includegraphics[width=0.75\textwidth]{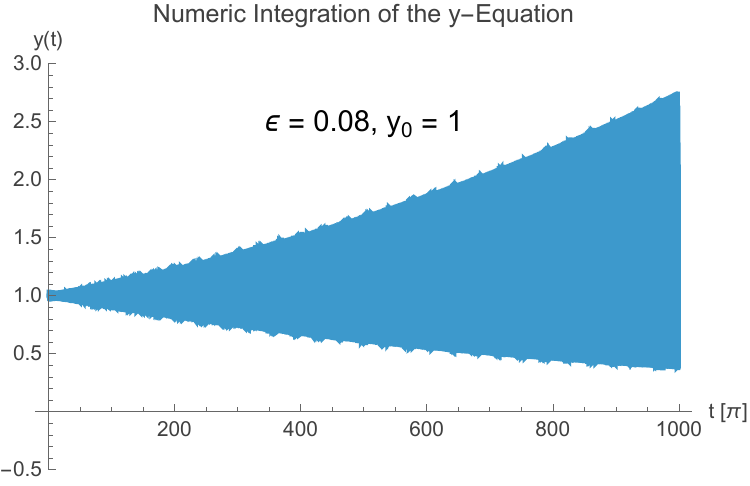}
	
	\vspace{0.8em}
	
	\includegraphics[width=0.75\textwidth]{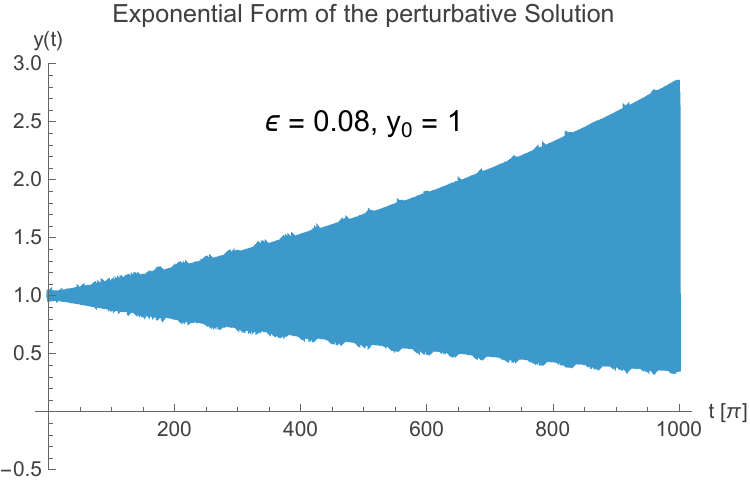}
	
	\caption{Comparison of the numerical solution of the $y$-equation (top)
		with the exponential perturbative approximation (bottom) for
		$\varepsilon=0.08$ and $y_0=1$.}
	\label{fig:y_num_exp}
\end{figure}

\begin{figure}[!htbp]
	\centering
	\includegraphics[width=0.75\textwidth]{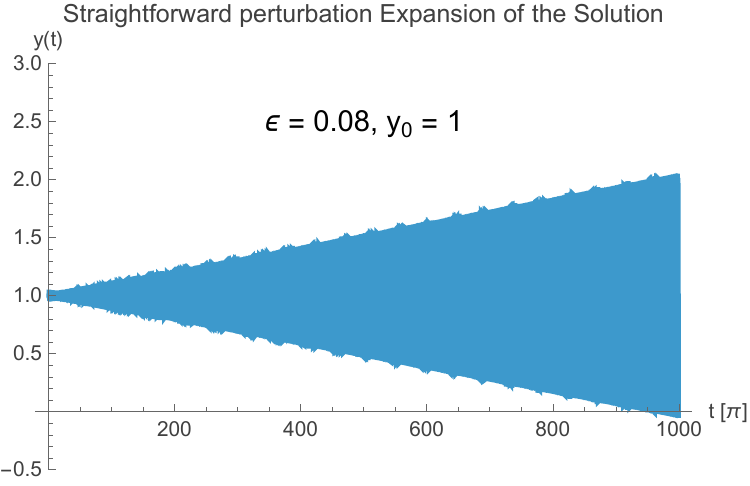}
	\caption{Straightforward perturbative expansion of the $y$-equation up to second order
		in $\varepsilon$. While accurate at early times, this approximation deteriorates at large
		times and does not preserve positivity.}
	\label{fig:y_taylor}
\end{figure}

Following (S1)-(S4) and inserting the perturbative
result for $y(t)$ into the general invariant framework yields an invariant for the solution of the z-equation of the form ($p=z'$):
\begin{equation}
	I(z,p,t)
	\;=\;
	A_1(t)\,z
	+ A_2(t)\,p
	+ A_3(t)\,z^2
	+ A_4(t)\,z p
	+ A_5(t)\,p^2
	+ A_6(t)\,z^3 -I(z_0,0,0)=0,
	\label{eq:I_general_form}
\end{equation}
where the coefficients $A_1$ to $A_6$ are given by 
\begin{equation}
	A_1(t)=\frac{\varepsilon}{2}\sin{t}   \label{eq:A1}
\end{equation}
\begin{equation}
	A_2(t)=\frac{\varepsilon}{2}\cos{t}   \label{eq:A2}
\end{equation}
\begin{equation}
	A_3(t)=y(t)+\frac{1}{2} y''(t)                \label{eq:A3}
\end{equation}
\begin{equation}
	A_4(t)=-y'(t)                         \label{eq:A4}
\end{equation}
\begin{equation}
	A_5(t)=y(t)                           \label{eq:A5}
\end{equation}
\begin{equation}
	A_6(t)=\frac{2}{3} y(t) g(t)                  \label{eq:A6}
\end{equation}

As can be seen the invariant $I(z,p,t)=K$ with $K=I(z_0,0,0)$ for any value of $t$ represents a quadratic equation in the momentum variable p as:
\begin{equation}
	A_5(t)p^2+(A_2(t)+A_4(t)z)p+A_6(t)z^3+A_3(t)z^2+A_1(t)z-K=0  \label{eq:p_quad}
\end{equation}
The usual two branched solution for $p_{1,2}(z,t)$ is then given by
\begin{equation}
p_{1,2}(z,t)=\frac{1}{2A_5} \left[-(A_2+A_4z)\pm \sqrt{A_2+A_4z)^2-4A_5(A_6z^3+A_3z^2+A_1z-K)} \right]  \label{eq:p12}
\end{equation}
Since the coefficients $A_1(t)$ - $A_6(t)$ contain secular contributions there is no time periodicity and the geometrical object related to the invariant becomes a tube with non periodic cross section. For $y(t)$ we inserted the second order perturbation expansion given by (\ref{eq:ypert}).  
\begin{figure}[htbp]
	\centering
	\begin{subfigure}[t]{0.48\textwidth}
		\centering
		\includegraphics[width=\linewidth]{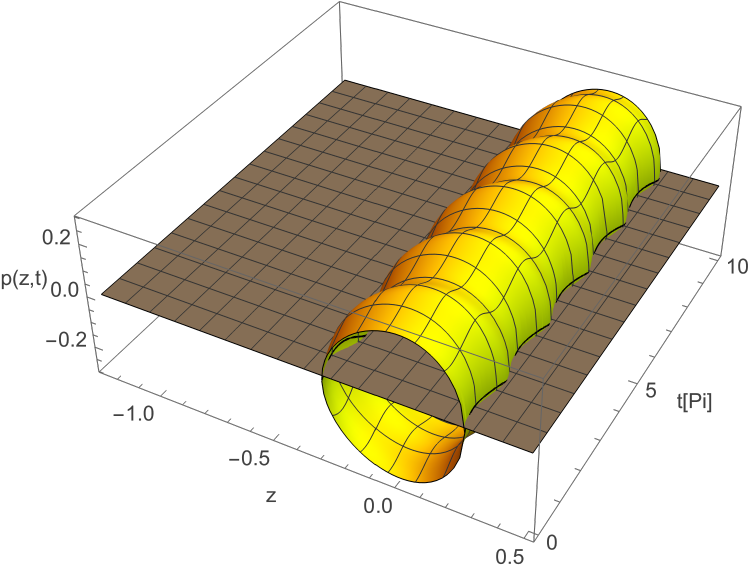}
		\caption{Approximate invariant surface $p=p(z,t)$ constructed from the perturbative approximation of $y(t)$ at early times ($0 \le t \le 10\pi$), showing a single-valued, closed tube geometry.}
		
		\label{fig:4a}
	\end{subfigure}\hfill
	\begin{subfigure}[t]{0.48\textwidth}
		\centering
		\includegraphics[width=\linewidth]{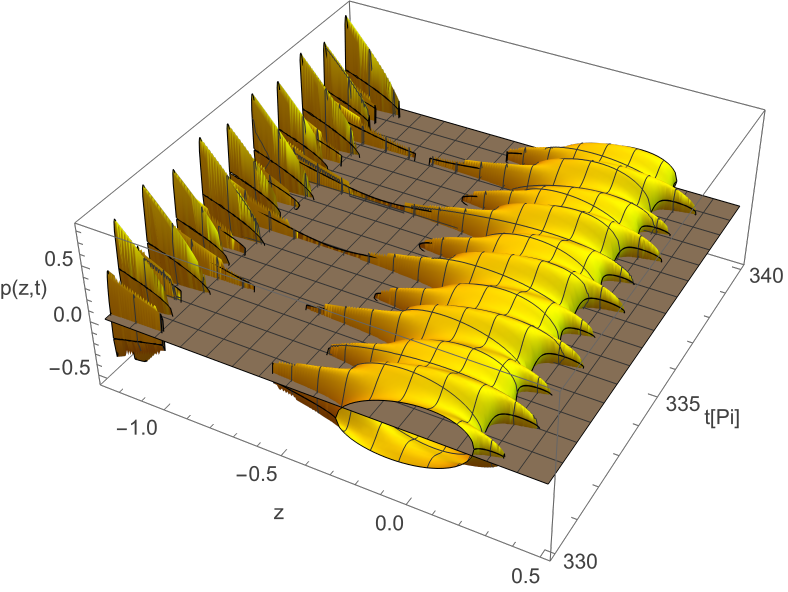}
		\caption{Same surface as in Fig.~4(a), shown at later times ($330\pi \le t \le 340\pi$), illustrating the onset of geometric overlap between neighboring branches.}

		\label{fig:4b}
	\end{subfigure}
	
	\caption{Geometry of the approximate invariant surface $p=p(z,t)$ at early and late times for $y_0=1$, $\varepsilon=0.08$ , $z_0=0.25$. While the surface remains single-valued and tube-like at early times, neighboring branches begin to overlap at later times, marking the geometric onset of the loss of confinement.
	}
	
	\label{fig:4}
\end{figure}
In Fig. 3 we demonstrate the approximate invariant surface for the parameters $y_0=1$ , $\varepsilon=0.08$ , $z_0=0.25$. We choose two cases, one for low values of time in the interval $0<t<10\pi$ and one for the advanced interval $320 \pi <t< 340\pi$. While in the first interval we see a well defined cylinder like surface with a nearly periodic cross section along the time axes, in the second interval the behaviour becomes more sophisticated. The central part of the invariant (tube), where the motion started at $t=0$ is still visible in a deformed way while a second part arising from the negative z-axes approaches and for short intevals even joins with the central tube. 
The depicted surface is constructed from the perturbative approximation of $y(t)$ and thus does not constitute an exact invariant of the full nonautonomous system.
Nevertheless, its geometry already reveals the mechanisms leading to the loss of confinement, which are subsequently quantified by the discriminant analysis. 
It should be noted that the occurrence of an overlap of the central tube with a structure reaching towards infinity for a finite, short period of time as shown in the figure, does not automatically imply an escape of the solution during this very time interval. This is because the actual orbit at the given time window could be located at a position on the surface not including the opening. At this point we may just state that the opening is a necessary condition for the occurrence of escape but not a sufficient one.

\subsection{Discriminant formalism and geometric classification}

Our analysis of rupture phenomena is based on the explicit algebraic structure of the approximate invariant introduced in the previous section.
Starting from the quadratic relation in the momentum variable $p$,
\begin{equation}
	A_5(t)p^2 + \bigl(A_2(t)+A_4(t)z\bigr)p
	+ A_6(t)z^3 + A_3(t)z^2 + A_1(t)z - K = 0,
\end{equation}
the existence of real solutions $p(z,t)$ is governed by the discriminant of this quadratic equation.
We denote this discriminant by
\begin{equation}
	\Disc(z,t).
\end{equation}

In the present work, $\Disc(z,t)$ is obtained by inserting the second--order exponential perturbation approximation for $y(t)$ into the invariant coefficients.
This representation is essential, since the exponential form preserves the positivity of $y(t)$ by construction and guarantees that $\Disc(z,t)$ remains a real--valued function for all times considered.

For every fixed time $t$, the function $z \mapsto \Disc(z,t)$ is a \emph{cubic polynomial in $z$}.
This observation is central: all geometric information about the opening or closure of the invariant tube at time $t$ is encoded in the real root structure of this cubic.

\paragraph{Algebraic characterization of bridge windows.}
The invariant surface defined by $I(z,p,t)=K$ is geometrically accessible only in those regions of $(z,p)$--space where $\Disc(z,t)\ge 0$.
Consequently, the topology of the set
\[
\{ z \in \mathbb{R} \mid \Disc(z,t)\ge 0 \}
\]
determines whether the invariant tube is closed or whether an unbounded passage in the $z$--direction exists.

This leads to the following purely algebraic classification.

\medskip
\noindent\textbf{Definition (Bridge window).}
A time $t$ belongs to a bridge window if the cubic equation
\begin{equation}
	\Disc(z,t)=0
\end{equation}
has \emph{exactly one real root}.

\medskip
The justification of this definition is straightforward and fundamental.
If the cubic possesses three real roots, the real axis is partitioned into multiple disjoint intervals separated by additional algebraic barriers, preventing a continuous connection to $|z|\to\infty$.
In contrast, if exactly one real root exists, the real line contains a connected half--line on which the sign of $\Disc(z,t)$ does not change.
In this case, a topologically open channel to infinity is present, and the invariant tube is locally open in the $z$--direction.

Hence, the condition ``exactly one real root'' is \emph{necessary and sufficient} for the existence of a bridge.

\paragraph{Reduction to a cubic discriminant criterion.}
The number of real roots of a cubic polynomial is completely determined by its cubic discriminant.
Denoting by
\[
\Delta_z\!\big[\Disc(z,t)\big]
\]
the discriminant of the cubic polynomial $\Disc(z,t)$ with respect to the variable $z$, the bridge condition can be expressed as the one--line criterion
\begin{equation}
	\label{eq:bridgecriterion}
	\mathrm{bridgeQ}(t)
	\;\Longleftrightarrow\;
	\Delta_z\!\big[\Disc(z,t)\big] < 0,
\end{equation}
up to the chosen sign convention.

While the quadratic discriminant is ubiquitous in elementary mathematics, the discriminant of a cubic polynomial rarely appears in dynamical applications.
In the present context, however, it provides a remarkably efficient and robust diagnostic tool: instead of analyzing the full two--dimensional discriminant surface $\Disc(z,t)$, one merely evaluates a scalar function of time.

\begin{remark}[Explicit form of the cubic discriminant]
	For completeness, we recall that the discriminant of a cubic polynomial
	\[
	P(z)=az^3+bz^2+cz+d
	\]
	is given explicitly by
	\begin{equation}
		\Delta(P)
		=
		18abcd - 4b^3d + b^2c^2 - 4ac^3 - 27a^2d^2.
	\end{equation}
	In the present setting, the coefficients $a,b,c,d$ are explicit time--dependent functions obtained from the quadratic discriminant $\Disc(z,t)$.
	Although these coefficients are known analytically, their full expansion is neither illuminating nor required for the analysis; only the sign of $\Delta_z\!\big[\Disc(z,t)\big]$ enters the bridge criterion \eqref{eq:bridgecriterion}.
\end{remark}

\paragraph{Geometric interpretation and relation to the exact system.}
Both the exact Bouquet--Hagel system and its perturbatively truncated counterpart exhibit rupture phenomena, but the underlying geometric mechanisms differ \cite{HagelBouquet1992,Hagel2025TubeRupt1}.

In the exact system, rupture occurs only for sufficiently large initial amplitudes $z_0$ and is associated with a gear--like interlocking of discriminant tongues.
In contrast, the truncated invariant studied here develops local bridge openings through the merging of opposing tongues at their tips, giving rise to isolated time intervals characterized by a single--root configuration of the cubic $\Disc(z,t)$.

The discriminant criterion \eqref{eq:bridgecriterion} isolates precisely this latter mechanism.
It does not contradict the exact theory; rather, it complements it by revealing an additional algebraic pathway to loss of confinement that is specific to perturbative truncation.
At this stage, the detailed shape of the discriminant surface and the mutual approach of its tongues are regarded as algebraic consequences of the time--dependent coefficients and will not be interpreted further.

\paragraph{Results.}
Applying the cubic discriminant criterion over extended time intervals yields a sequence of bridge windows that can be represented compactly as a box plot along the time axis.
When time is measured in units of $\pi$, the resulting structure displays a clear organization: initially, bridge windows are separated by intervals of length $2\pi$, corresponding to a single active phase per cycle.
Beyond a critical time $t/\pi \approx 372$, this separation collapses to $\pi$, indicating a period--doubling in the bridge timing.
Moreover, the widths of successive windows alternate, reflecting the emergence of two inequivalent phase classes within one $2\pi$--cycle.

A direct comparison with the discriminant surface $\Disc(z,t)$ over the same time interval confirms that the box plot provides a compressed but faithful representation of the underlying geometry.
While plotting $\Disc(z,t)$ offers valuable visual intuition, it is computationally expensive and conceptually unnecessary: the rupture mechanism is entirely captured by the cubic discriminant criterion.

\begin{figure}[t]
	\centering
	\includegraphics[width=\linewidth]{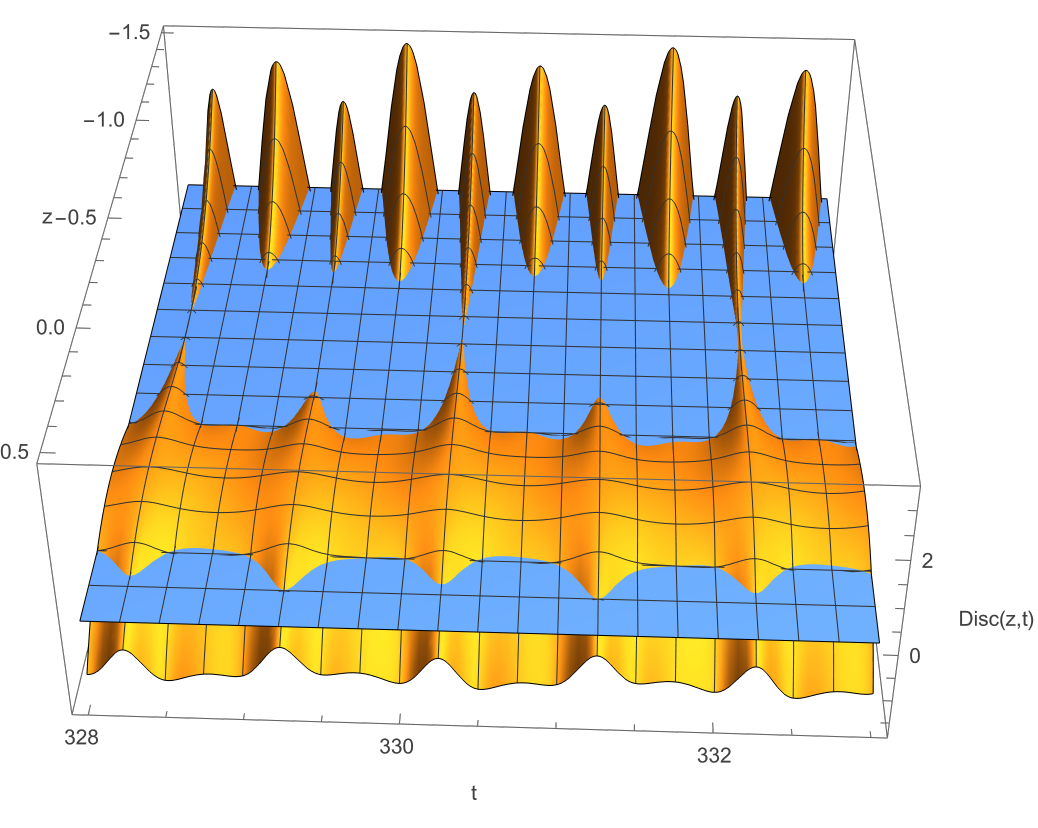}
	\caption{Onset of the bridge phenomenon in the discriminant surface $\Disc(z,t)$. 
		For $t\approx 328\pi$ no bridge is present and the invariant tube remains closed. 
		Around $t\approx 330\pi$ the discriminant tongues nearly touch, indicating the imminent opening of a bridge. 
		At $t\approx 332\pi$ a clear local opening is visible, corresponding to the first time interval in which the cubic equation $\Disc(z,t)=0$ possesses exactly one real root.}
	\label{fig:bridge_onset}
\end{figure}

\begin{figure}[t]
	\centering
	\begin{subfigure}[t]{0.48\linewidth}
		\centering
		\includegraphics[width=\linewidth]{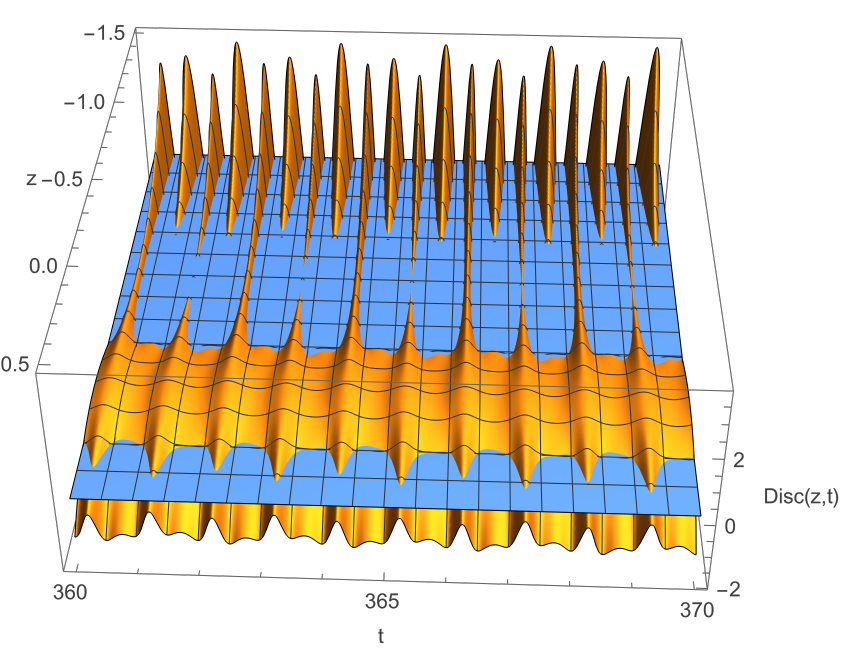}
		\caption{Discriminant surface $\Disc(z,t)$ in the interval $360\le t/\pi\le370$. 
			A single bridge window per $2\pi$ cycle is still present, but the discriminant tongues begin to deform and flatten, marking the onset of a period--doubling transition.}
		\label{fig:doubling_start}
	\end{subfigure}
	\hfill
	\begin{subfigure}[t]{0.48\linewidth}
		\centering
		\includegraphics[width=\linewidth]{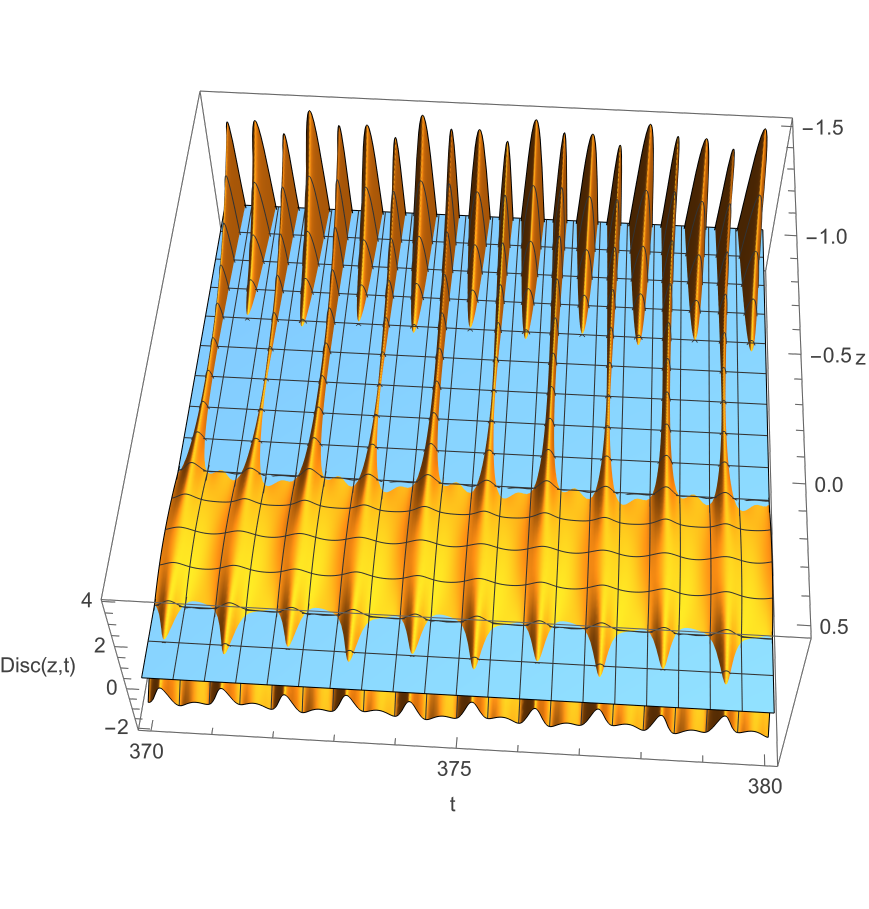}
		\caption{Same representation for $370\le t/\pi\le380$. 
			Two distinct bridge windows per $2\pi$ cycle are now clearly visible, corresponding to a separation of $\pi$ in time. The alternating shape and width of the openings reflect the emergence of two inequivalent phase classes.}
		\label{fig:doubling_full}
	\end{subfigure}
	
	\caption{Emergence of period--doubling in the bridge geometry of the discriminant surface $\Disc(z,t)$. 
		Panel (a) shows the last regime with a single bridge per $2\pi$ cycle, while panel (b) demonstrates the fully developed period--halved structure. 
		Careful inspection reveals the formation of secondary bulges in the central channel, indicating early geometric signatures of a subsequent doubling, although no further transition is claimed at this stage.}
	\label{fig:period_doubling}
\end{figure}

While Figs. 4 and 5 visualize the geometric opening of bridge windows in the discriminant surface, this information can be condensed further.
Since the bridge criterion is purely algebraic and depends only on time, it allows for a compact representation of all escape–enabling intervals along the time axis.
In the following, we therefore introduce a box plot that encodes the occurrence and duration of bridge windows without reference to the full $(z,t)$–geometry.

\begin{figure}[htbp]
	\centering
	\includegraphics[width=\linewidth]{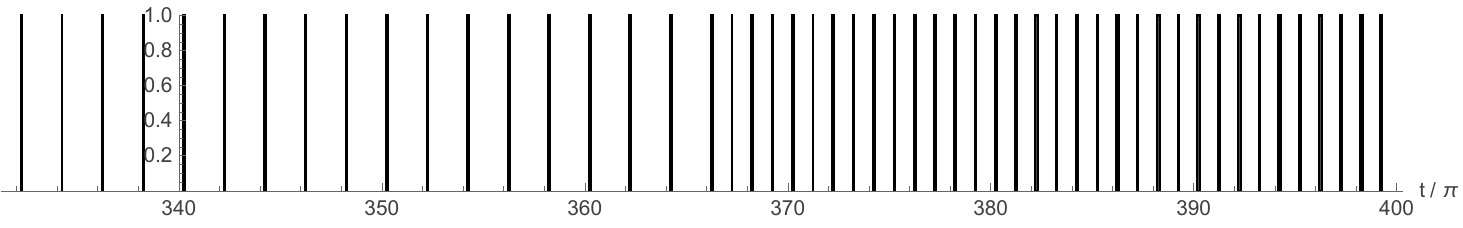}
	\caption{Box plot representation of bridge windows obtained from the cubic discriminant criterion.
		Each vertical box marks a time interval in which the cubic equation $\Disc(z,t)=0$ possesses exactly one real root.
		The horizontal axis shows time in units of $t/\pi$, while the vertical axis represents the Boolean bridge indicator (taking values $0$ or $1$), included only for visual separation.
		Around $t/\pi\approx369$ the spacing between successive boxes changes from $2\pi$ to $\pi$, signaling the onset of period--doubling in the bridge timing.
		At the beginning of this transition the newly appearing $\pi$--spaced boxes are still very narrow, but they subsequently broaden in a manner similar to the original $2\pi$--periodic windows, indicating an increasing temporal persistence of individual bridge intervals at later times.}
	\label{fig:boxplot}
\end{figure}

Although no probabilistic model is invoked, the increasing visual density of the box plot illustrates how escape-enabling intervals occupy a progressively larger fraction of the time axis.

The box plot provides a compact time–resolved view of all bridge windows detected by the cubic discriminant criterion.
It clearly reproduces the period–doubling observed in the discriminant surfaces: the separation between successive windows changes abruptly from $2\pi$ to $\pi$.
In addition, the widths of the boxes exhibit a gradual increase with time, indicating that individual bridge windows persist longer at later stages.

\FloatBarrier

\section{Conclusions}

The original motivation of this work was to obtain a controlled approximation of an
exactly tube-integrable nonautonomous system and to assess how faithfully its geometric
structure can be reproduced beyond the strictly integrable setting \cite{Hagel2025TubeInteg}.

To this end, we constructed a perturbative approximation of the auxiliary $y$--equation . 
governing the time dependence of the invariant coefficients, enforcing positivity by means
of an exponential representation.
This approach was designed to preserve the geometric admissibility of the invariant,
while systematically improving the accuracy of the approximation compared to earlier
treatments.

Within this refined framework, an unexpected phenomenon emerges.
Rather than a gradual loss of confinement, the approximate invariant develops localized
rupture intervals in time, characterized by the appearance of algebraic bridges in the
invariant tube \cite{Hagel2025TubeRupt1}.

These bridge windows exhibit a nontrivial temporal organization, including a transition
from one opening per $2\pi$ cycle to two openings per cycle, reminiscent of a period-doubling
scenario.
Although no probabilistic interpretation is implied, the resulting window structure shows
formal similarities to Bernoulli-type intermittency, motivating the terminology adopted
in this work. The terminology is used in a purely descriptive sense and does not imply
a dynamical equivalence with classical intermittency scenarios
\cite{Berry1989}.

A central result is that this rupture behavior can be detected and classified by a purely
algebraic criterion.
By reducing the geometric problem to the sign of the discriminant of a cubic polynomial,
the onset and temporal evolution of tube opening are captured by a single scalar function
of time.
This reduction remains effective despite the weak non-integrability introduced by
perturbative truncation and highlights the robustness of algebraic diagnostics derived
from time-dependent invariants.

An important methodological distinction from earlier studies concerns the role of
sampling.
In previous work on weakly non-integrable approximations, only the secular contributions
to the invariant coefficients were retained.
In that reduced setting, no phase drift was present, and sampling at integer multiples of
the base period intersected the discriminant tongues in a controlled and analytically
tractable manner.
In the present, more accurate approximation, additional oscillatory terms are retained.
As a consequence, a slow phase shift is implicitly introduced into the dynamics of the
auxiliary equation.
Although this phase drift is not yet isolated in explicit form, its presence prevents a
direct extension of the earlier sampling strategy.

Rather than a drawback, this limitation points toward a continuation of the
analysis.
A systematic extraction of the phase drift from the perturbative solution of the
$y$--equation would allow the sampling method to be reinstated in a modified form and may
lead to particularly simple analytic expressions for rupture thresholds and bridge window
boundaries.
More generally, the present results indicate that increasingly accurate approximations of
time-dependent invariants can reveal qualitatively new geometric phenomena, such as the
emergence and densification of rupture windows, even when exact tube integrability is no
longer preserved.

The present results further support the view of the Bouquet--Hagel system as a versatile
reference framework, whose free parameters allow meaningful adaptations to a range of
time-dependent problems, as already demonstrated in earlier applications to celestial
mechanics.

\section*{Acknowledgements}

The author gratefully acknowledges long-standing discussions with C.~Bouquet,
which originally led to the development of the integrable reference system
considered in this work and to its later applications.
Parts of the present analysis build upon methods developed in earlier joint and
individual contributions.

The author also acknowledges the use of the AI language model ChatGPT (version~5.2)
as a supportive tool for technical editing, consistency checks, and the refinement
of mathematical exposition.
All scientific content, interpretations, and conclusions are solely the
responsibility of the author.

\end{document}